\documentstyle{amsart}

\def\P{{\hbox{\bf P}}}
\def\E{{\hbox{\bf E}}}

 at 10 true pt

\def\be#1{ \begin{equation}\label{#1} }
\def\bas{\begin{align*}}
\def\eas{\end{align*}}
\def\bi{\begin{itemize}}
\def\ei{\end{itemize}}

\newenvironment{proof}{\noindent {\bf Proof} }{\endprf\par}
\def \endprf{\hfill  {\vrule height6pt width6pt depth0pt}\medskip}
\def\emph#1{{\it #1}}
\def\textbf#1{{\bf #1}}

\def\CQ{{\mathcal Q}}

\def\
tl{{\tilde l}}

\def\ep{{\epsilon}}

\theoremstyle{plain}

\newtheorem{theorem}[subsection]{Theorem}

\newtheorem{fact}[subsection]{Fact}
\newtheorem{lemma}[subsection]{Lemma}

\newtheorem{corollary}[subsection]{Corollary}

\newtheorem{question}[subsection]{Question}
\theoremstyle{remark}

\newtheorem{remark}[subsection]{Remark}

\theoremstyle{definition}
\newtheorem{definition}[subsection]{Definition}

\begin{document}

\title{The Rank of Random Graphs}
\author{Kevin P. Costello}

\address{Department of Mathematics, UCSD, La Jolla, CA 92093}
\email{kcostell@@ucsd.edu}
\thanks{Kevin P. Costello is supported by
an NSF Graduate Research Fellowship}

\author{Van H. Vu}
\address{Department of Mathematics, Rutgers, Piscataway, NJ 08854}

\email{vanvu@@math.rutgers.edu}
\thanks{Van H. Vu is a Sloan Fellow and is supported by an NSF Career Grant.}

\maketitle
\begin{abstract}
We show that almost surely the rank of the adjacency matrix of the
Erd\"os-R\'enyi random graph $G(n,p) $ equals the number of
non-isolated vertices for any $c\ln n/n<p<1/2$, where $c$ is an
arbitrary positive constant larger than $1/2$.
 In particular the giant component (a.s.) has full rank in this range.
\end{abstract}

\section{Introduction}

Let $G$ be a (simple) graph on $n$ points $\{1, \dots, n\}$. The
adjacency matrix $Q_G$ of $G$ is a symmetric $n$ by $n$ matrix,
whose $ij$ entry is one if the vertex $i$ is connected to the
vertex $j$ and zero otherwise.

The rank of $Q_G$ is a basic  graph parameter of the graph and has
been used  in graph theory and also in applications from computer
science. Thus, it is natural and perhaps important to understand
the behavior of this parameter with respect to a random graph.

An interesting feature of this parameter is that, unlike most of
graph parameters (such as the connectivity, the chromatic number
etc), the rank is not monotone (adding edges may reduce the rank).
So, it is not clear that one can define thresholds. Our study,
however,  will shed some light on this problem.

There are several models for random graphs. We will focus on the
most popular one, the Erd\" os-R\'enyi $G(n,p)$ model (some other
models will be discussed in the concluding remarks). In this model,
one starts with the vertex set $\{1, \dots, n\}$ and puts an edge
(randomly and independently) between any two distinct vertices $i$
and $j$ with probability $p$. We say that a property $P$ holds
almost surely for $ G(n,p)$ if the probability that $G(n,p)$
possesses $P$ goes to one as $n$ tends to infinity.

A vertex $v$ of a graph $G$ is isolated if it has no neighbor. If
$v$ is isolated, then the row corresponding to $v$ in $Q_G$ is
all-zero. Let $i(G)$ denote the number of isolated vertices of
$G$.  It is clear that

\begin{fact} \label{fact:1} For any graph $G$, ${\hbox {rank}} (Q_G) \le n -i(G). $
\end{fact}

The main result of this paper shows that for $p$ sufficiently
large, the above upper  bound is tight for $G(n,p)$. In other
words, any non-trivial linear dependence in $G(n,p)$ comes from
the isolated vertices.

\begin{theorem} \label{theo:main}
Let $c$ be a constant larger than $\frac{1}{2}$.  Then for any $
 c \ln n / n \le p \le \frac{1}{2}$, the following holds with probability $1- O((\ln
\ln n)^{-1/4})$  for a random sample $G$ from $G(n,p)$:

$${\hbox {rank}} (G) = n - i(G). $$
\end{theorem}

\begin{remark}
This result is sharp in two ways. Obviously, the estimate $\\ {\hbox
{rank}} (G) =n -i(G)$ cannot be improved. Furthermore, the bound $c
> 1/2$ is also the best possible. For $c < 1/2$ and $p = c \ln n/n$,
a random sample $G$ from $G(n,p)$ satisfies the strict inequality

$${\hbox {rank}} (G) < n- i(G) $$

\noindent almost surely. In order to see this, notice that in this
range $G(n,p)$ almost surely contains two vertices $u$ and $v$ with
degree one sharing a common neighbor. The rows corresponding to $u$
and $v$ are not zero, but they are equal and this reduces the rank
further.
\end{remark}

Let us now deduce a few corollaries. It is well known that $p \ge
c\ln n/n$ for some $c >1/2$, then the random graph (a.s.) consists
of a giant component and some isolated vertices.

\begin{corollary} \label{cor:main1} Let $c$ be a constant larger than $\frac{1}{2}$. Then for
any $c \ln n / n< p<1/2$, the giant component of $G(n,p)$ is almost
surely non-singular.
\end{corollary}

\noindent Furthermore,  if $c > 1$, then $G(n,p)$ almost surely is
connected and contains no isolated vertices.

\begin{corollary} \label{cor:main} Let $c$ be a constant larger than 1. Then for any
$c \ln n / n< p<1/2$, $G(n,p)$ is almost surely non-singular.
\end{corollary}

It follows that the non-singularity of $G(n,p)$ has a sharp
threshold at $p =\ln n/n$. For any positive constant $\ep$, if $p
< (1-\ep) \ln n/n$, then $G(n,p)$ is almost surely singular as it
contains isolated vertices. On the other hand, the above corollary
asserts that $G(n,p)$ is almost surely non-singular for $p >
(1+\ep) \ln n/n$.

The  special case $p=1/2$ was a well known  conjecture of B.
Weiss, posed many years ago. This special case can be viewed as
the symmetric version of a well-known theorem of Koml\'os on the
non-singularity of (non-symmetric) random Bernoulli matrices
\cite{Kom1}  and was solved two years ago in \cite{CTV}. The proof
of Theorem \ref{theo:main} extends the ideas in that paper and
combines them with arguments involving random graphs. In the next
section, we outline this proof and present the key lemmas. In
Section 3, we prove the theorem assuming the lemmas. Most of the
rest of the paper is devoted to the proofs of the lemmas. The last
section contains several remarks and open questions.

{\bf Notation.} In the whole paper, we will always assume that $n$
is sufficiently large. As usual, the asymptotic notation is used
under the condition that $n \rightarrow \infty$. $\P$ denotes
probability and $\E$ denotes expectation.

\section{Outline of the proof and the main lemmas}

We assume that $o(1) = p \ge c\ln n/n$ for a constant $c >1/2$. The
treatment of the case $p=\Omega (1)$  was presented in \cite{CTV}
and is  omitted here. We denote by $Q(n,p)$ the adjacency matrix of
$G(n,p)$.

Following the ideas from \cite{Kom1, CTV}, we are going to expose
$Q(n,p)$
minor by minor. Letting $Q_m$ denote the upper left $m \times m$ minor of $%
Q(n,p)$, we view $Q_{m+1}$ as being formed by taking $Q_m$ and
augmenting by a column whose entries are chosen independently,
along with the column's transpose.  Denote by $G_m$  the graph
whose adjacency matrix is $Q_m$.  In graph theoretic terms, we are
considering the vertex exposure process of $G(n,p)$.

Our starting observation is that when a good portion of the vertices have been
exposed, the matrix has rank close to its size.

Recall that $p \ge c\ln n/ n$ for a constant $c >1/2$. We can set
a constant $0 <\delta<1$ such that $1/2<\delta c < 3/5$. Define
$n^{\prime}:=\delta n$.

\begin{lemma} \label{lemma:nearfullrank} For any constant $\epsilon>0$ there
exists a constant $\gamma>0$ such that
\begin{equation*}
{\hbox{\bf P}}(\hbox{rank}(Q_{n^{\prime}})<(1-\epsilon)n^{\prime})=o(e^{-%
\gamma n \ln n})
\end{equation*}

\end{lemma}

Our plan is to show that the addition of the remaining
$n-n^{\prime}$
rows/columns is enough to remove all the linear dependencies  from $%
Q_{n^{\prime}}$, except those corresponding to the isolated
vertices.

The next lemmas provide some properties of $G_m$ for
$n^{\prime}\leq m \leq n$.

\begin{definition}
A graph $G$ is {\bf well-separated} if it contains no pair of
vertices of degree at most $\ln \ln n$ whose distance from each
other is at most 2.
\end{definition}

\begin{lemma} \label{lemma:separation} For any constant $\ep >0$, $G_m$ is well
separated for every $m$ between $n'$ and $n$ with probability
$1-O(n^{1-2c\delta+\ep})$.
\end{lemma}

Here and later on, we always choose $\ep$ sufficiently small so
that $1 -2c\delta +\ep$ is negative.

\begin{definition}

A graph $G$ is a {\bf small set expander} if every subset $S$ of the
vertices of $G$ with $|S| \leq \frac{n}{\ln^{3/2} n}$ containing no
isolated vertices has at least $|S|$ edges connecting $S$ to
$\bar{S}$, its complement.
\end{definition}

\begin{lemma} \label{lemma:expansion}

For any $m>n'$ the probability that $G_m$ is well separated but is
not a small set expander is $O(1/n^3)$.
\end{lemma}

\begin{remark}
Lemmas \ref{lemma:separation} and \ref{lemma:expansion} immediately
imply that almost every graph encountered in our process after time
$n'$ will be a small set expander. We cannot expect this to occur
for $p<\frac{(.5-\ep)\ln n}{n}$, as at this density  the random
graph will likely contain pairs of adjacent vertices of degree 1,
leading to sets of size 2 without any edges at all leaving them.
\end{remark}

\begin{definition} \label{def:distinguished} A set $S$ of the vertices of a graph $G$ is {\bf %
nice} if there are at least two vertices of $G$ each is adjacent
to exactly one vertex in $S$.
\end{definition}

Set $k := {\frac{\ln \ln n }{2 p}}$.

\begin{definition} \label{def:goodgraph} A graph $G$ is {\bf good} if the following two
properties hold:

1. Every subset of the vertices of $G$ of size at least 2 and at
most $k$ which contains no isolated vertices is nice.

2. At most $\frac{1}{p \ln n}$ vertices of $G$ have degree less than
2.

A symmetric (0,1) matrix $A$ is {\bf good} if the graph for which it
is an adjacency matrix is good.
\end{definition}

The next lemma states that in the  augmentation process we will
likely run only into good matrices.
\begin{lemma}

\label{lemma:goodmatrices} Let $\ep$ be a positive constant. Then
with probability $1-O(n^{1-2c\delta + \epsilon})$, $Q_m$ is good
for every $m$ between $n'$ and $n$ .
\end{lemma}

Finally, we have two lemmas stating that good matrices behave nicely
when augmented.

\begin{definition} \label{def:normalpair}

A pair $(A, A')$ of matrices is called {\bf normal} if $A'$ is an
augmentation of $A$ and every row of all 0's in $A$ also contains
only 0's in $A'$ (in graph theoretic terms, the new vertex added by
the augmentation is not adjacent to any vertices which were isolated
before the augmentation).
\end{definition}

\begin{lemma} \label{lemma:singaug} Let $A$ be any fixed, good $m \times m$ matrix
with the property that $rank(A)+i(A)<m$. Then
\begin{equation*}
{\hbox{\bf P}}(\,\, \hbox{rank}(Q_{m+1})-\,\, \hbox{rank}(Q_m)<2 |
(Q_m, Q_{m+1}) \, \, \hbox{is normal} \wedge
Q_m=A)=O((kp)^{-1/2}).
\end{equation*}

\end{lemma}

\begin{lemma} \label{lemma:nonsingaug} Let $A$ be any fixed, good, $m \times m $
matrix with the property that $rank(A)+i(A)=m$. Then
\begin{equation*}
{\hbox{\bf P}}(\,\, \hbox{rank}(Q_{m+1})-\,\, \hbox{rank}(Q_m)<1 |
(Q_m, Q_{m+1}) \, \, \hbox{is normal} \wedge
Q_m=A)=O((kp)^{-1/4}).
\end{equation*}

\end{lemma}

Set  $Y_m= m- (rank(Q_m)+i(Q_m))$. The above two lemmas force
$Y_m$ to stay near 0. Indeed, if $Y_m$ is positive, then when the
matrix is augmented, $m$ increases by 1 but the rank of $Q_m$ will
likely increase by 2 (notice that $kp \rightarrow \infty$),
reducing $Y_m$. On the other hand, if $Y_m =0$,  it is likely to
stay the same after the augmentation.

 In the next section, we will turn this heuristic into a rigorous calculation and prove  Theorem
 \ref{theo:main}, assuming the  lemmas.

\section{Proof of the Main Result from the Lemmas}

In this section, we  assume all lemmas are true. We are going to
use a variant of an argument from \cite{CTV}.

Let $B_1$ be the event that the rank of $Q_{n^{\prime}}$ is at
least $n^{\prime}(1-{\frac{1-\delta }{4 \delta}})$.  Let $B_2$ be
the event that $Q_m$ is good for all $n^{\prime}\leq m < n$.  By
Bayes' theorem we have
\begin{equation*}
{\hbox{\bf P}}(\,\, \hbox{rank}(Q_n)+i(Q_n)<n) \leq {\hbox{\bf P}}(\,\, \hbox{rank}%
(Q_n)+i(Q_n)<n \wedge B_2 | B_1) + {\hbox{\bf P}}(\neg B_1) +
{\hbox{\bf P}}(\neg B_2)
\end{equation*}

By Lemma \ref{lemma:nearfullrank} we have that ${\hbox{\bf P}}(\neg
B_1)=o(e^{-\gamma n \ln n})$ and by Lemma \ref{lemma:goodmatrices}
${\hbox{\hbox{\bf P}}(\neg B_2)=O(n^{1-2c\delta +\ep})}$. Both
probabilities  are thus much smaller than the bound $O((\ln \ln
n)^{-1/4})$ which  we are trying to prove. So, it remains to bound
the first term.

Let $Y_m=m-\,\, \hbox{rank}(Q_m)-i(Q_m)$. Define a random variable
$X_m$ as follows:

\begin{itemize}

\item  $X_m=4^{Y_m}$ if $Y_m>0$ and every $Q_j$ with $n' \leq j
\leq m$ is good;

\item  $X_m=0$ otherwise.

\end{itemize}

The core of the proof is the following  bound on  the expectation
of $X_{m+1}$ given any fixed sequence  $\CQ_m$ of matrices
$\{Q_{n'}, Q_{n'+1}, \dots, Q_m\}$ encountered in the augmentation
process.

\begin{lemma} \label{lemma:expectationofX} For any sequence $ \CQ_m= \{Q_{n'}, Q_{n'+1}, \dots, Q_m\}$
encountered in the augmentation process,

\begin{equation*}
{\hbox{\bf E}} (X_{m+1} |  \CQ_m  ) \leq \frac{3}{5} X_m + O((\ln
\ln n)^{-1/4}).
\end{equation*}

\end{lemma}

Let us (for now) assume Lemma \ref{lemma:expectationofX} to be true.
This lemma together with Bayes theorem shows that for $n^{\prime}<m$
we have
\begin{equation*}
{\hbox{\bf E}}(X_{m+1} | Q_{n'})<{\frac{3 }{5}}{\hbox{\bf E}}(X_m
| Q_{n'})+O((\ln \ln n)^{-1/4}).
\end{equation*}

By induction on $m_2-m_1$ we now have that for any $m_2\geq m_1 \geq
n^{\prime}$
\begin{equation*}
{\hbox{\bf E}}(X_{m_2} | Q_{n'})<(\frac{3}{5})^{m_2-m_1} {\hbox{\bf E}}%
(X_{m_1} | Q_{n'})+O((\ln \ln n)^{-1/4}).
\end{equation*}

In particular, by taking $m_2=n$ and $m_1=n^{\prime}$ we get that
\begin{equation*}
{\hbox{\bf E}}(X_n | Q_{n'})<(\frac{3}{5})^{n-n'} X_{n^{\prime}} +
O((\ln \ln n)^{-1/4}).
\end{equation*}

\noindent If $Q_{n'}$ satisfies $B_1$, we automatically have
$X_{n'} \leq 4^{\frac{(1-\delta)n'}{4 \delta}}=(\sqrt{2})^{n-n'}$,
so
\begin{equation*}
{\hbox{\bf E}}(X_n | Q_{n'}) <(\frac{3 \sqrt{2}}{5})^{n-n'} +
O((\ln \ln n)^{-1/4})=O((\ln \ln n)^{-1/4}).
\end{equation*}

\noindent By Markov's inequality, for any $Q_{n'}$ satisfying
$B_1$
\begin{equation*}
{\hbox{\bf P}}(X_n>3 | Q_{n'})=O((\ln \ln n)^{-1/4})
\end{equation*}

\noindent On the other hand, by definition $X_n\ge 4$ if
${\hbox{rank}}(Q_n)+i(Q_n)<n$ and $B_2$ holds. It thus follows by
summing over all $Q_{n'}$ satisfying $B_1$ that

\begin{equation*}
{\hbox{\bf P}}(\,\, \hbox{rank}%
(Q_n)+i(Q_n)<n \wedge B_2 | B_1) =O((\ln \ln n)^{-1/4}),
\end{equation*}

\noindent proving the theorem.

\vskip2mm

It remains to prove Lemma \ref{lemma:expectationofX}. If a matrix
in the sequence $\CQ_m= \{Q_{n'}, Q_{n'+1}, \dots, Q_m\}$ is not
good, then $X_{m+1}=0$ by definition and there is nothing to
prove. Thus, from now on we can assume that all matrices in the
sequence are good.

Let $Z_m$ denote the number of vertices which were isolated in
$Q_m$ but not in $Q_{m+1}$.  If $Z_m$ is positive, then augmenting
the matrix will increase $Y_m$ by at most $Z_m+1$ ($m$ increases
by 1, the number of isolated vertices decreases by at most  $Z_m$,
and the rank does not decrease).  Furthermore, $Z_m=0$ if and only
if $(Q_m, Q_{m+1})$ is normal.

By Bayes' theorem, we have
\begin{eqnarray*}
{\hbox{\bf E}}(X_{m+1}| \CQ_m) &=& {\hbox{\bf E}} (X_{m+1}| Z_m >0
\wedge \CQ_m) {\hbox{\bf P}} (Z_m>0 |\CQ_m)
  \\&\,\,+&{\hbox{\bf E}}
(X_{m+1}| \CQ_m \wedge (Q_m, Q_{m+1}) \, \, \hbox{is normal}) {\hbox{\bf P}} ((Q_m, Q_{m+1}) \, \, \hbox{is normal}| \CQ_m)  \\
&\leq& {\hbox{\bf E}} (X_{m+1} \, \chi (Z_m>0)| \CQ_m) +{\hbox{\bf
E}} (X_{m+1}| \CQ_m \wedge (Q_m, Q_{m+1}) \, \, \hbox{is normal})
\\
&=&  {\hbox{\bf E}} (4^{Z_m +1+Y_m} \chi (Z_m>0) |\CQ_m) +
{\hbox{\bf E}} (X_{m+1}| \CQ_m \wedge (Q_m, Q_{m+1}) \, \,
\hbox{is normal}) .
\end{eqnarray*}

Since $Q_m$ is good, $G_m$ has at most $\frac{1}{p \ln n}$
isolated vertices. Thus,  we can bound $Z_m$ by the  sum of
$\frac{1}{p \ln n}$ random Bernoulli variables, each of which is 1
with probability $p$. It follows that
\begin{equation*}
\hbox{\bf P} (Z_m = i) \leq \binom{(p \ln n)^{-1}}{i} p^{i} \leq
(\ln n)^{-i}.
\end{equation*}

Adding up over all $i$, we have
\begin{equation*}
\hbox{\bf E} (4^{Z_m + 1} \chi(Z_m>0 )|\CQ_m) \leq
\sum_{i=1}^{\infty} 4^{i+1} (\ln n)^{-i} = O(({\ln n})^{-1}).
\end{equation*}

If $Y_m=0$ and $(Q_m, Q_{m+1})$ is normal, then  by Lemma
\ref{lemma:nonsingaug} (which applies since $Q_m$ is good)
$X_{m+1}$ is either 0 or 4, with the probability of the latter
being $O((\ln \ln n)^{-1/4})$. Therefore we have for any sequence
$\CQ_m= \{Q_{n'}, \dots Q_m\}$ of good matrices with $Y_m=0$ that
\begin{equation} \label{fullrankexpec}
{\hbox{\bf E}}(X_{m+1} | \CQ_m)=O((\ln \ln n)^{-1/4}+(\ln
n)^{-1})=O((\ln \ln n)^{-1/4}).
\end{equation}

If $Y_m=j>0$ and $(Q_m, Q_{m+1})$ is normal, then $Y_{m+1}$ is
$j-1$ with probability $\newline 1-O((\ln \ln n)^{-1/2})$ by Lemma
\ref{lemma:singaug}, and otherwise is at most $j+1 $.  Combining
this with the bound on $\hbox{\bf E} (4^{Z_m + 1} \chi(Z_m>0
)|\CQ_m)$ we have

\begin{equation} \label{nonfullrankexpec}
{\hbox{\bf E}} (X_{m+1}|\CQ_m)=4^{j-1}+4^{j+1}
O((\ln \ln n)^{-1/2})+4^j O((\ln n)^{-1}) \leq {\frac{3 }{5}%
}4^j
\end{equation}

The lemma now follows immediately from  (\ref{fullrankexpec}) and
(\ref{nonfullrankexpec}).

\section{Proof of Lemma \protect\ref{lemma:nearfullrank}}

By symmetry and the union bound

 $$\P( \text{rank%
}(Q_{n^{\prime}})<(1-\epsilon) n^{\prime}) \le
\binom{n^{\prime}}{\epsilon n^{\prime}} \times \P(B_1^*),$$ where
$B_1^*$ denotes the event  that the last $\epsilon n^{\prime}$
columns of $Q_n^{\prime}$ are contained in the span of the remaining
columns.

We view $Q_{n'}$ as a block matrix,
$$Q_{n'}=\left[ \begin{array}{c|c} A & B \\ \hline B^T & C \\ \end{array}
\right],$$ where $A$ is the upper left $(1-\ep)n' \times (1-\ep)n'$
sub-matrix and $C$ has dimension $\ep n' \times \ep n'$.  We obtain
an upper bound on $\P(B_1^*)$ by bounding the probability of $B_1^*$
conditioned on any fixed $A$ and $B$ (treating $C$ as random).

$B_1^*$ cannot hold unless the columns of $B$ are contained in the
span of those of $A$, meaning the equation $B=A F$ holds for some
matrix $F$.  If this is the case, then $B_1^*$ will hold only when
we also have $C=B^T F$.  This means that each entry of $C$ is forced
by our choice of $A$, $B$ and our assumption that $B_1^*$ holds.

However, $C$ is still random, and the probability that any given
entry takes on its forced value is at most $1-p$.  The entries are
not all independent (due to the symmetry of $C$), but those on or
above the main diagonal are.  Therefore the probability that $B_1^*$
holds for any fixed $A$ and $B$ is at most $(1-p)^{\frac{(\ep
n')^2}{2}}$.

We therefore have

\begin{eqnarray*}
{\hbox{\bf P}} (\text{rank}(Q_{n^{\prime}})<(1-\epsilon) n^{\prime})
&\leq& \binom{n^{\prime}}{\epsilon n^{\prime}}
((1-p)^{\frac{(\epsilon n^{\prime})^2 }{2}} \\
&\leq& ({\frac{n^{\prime}e }{\epsilon n^{\prime}}})^{\epsilon
n^{\prime}} e^{\frac{-
p (\epsilon n^{\prime})^2 }{2}} \\
&\leq& {c_2}^{n} e^{-c_1 n \ln n}.
\end{eqnarray*}

where $c_1$ and $c_2$ are positive constants depending on
$\epsilon$, $\delta $, and $c$ (but independent of $n$).

\begin{remark} \label{remark:smallerp} The same argument gives an upper bound of
${c_2}^n e^{-c_1 n^2 p}$ on the probability for any $n$ and $p$.
Holding $\epsilon$ fixed, we see that the probability becomes
$o(1)$ for $p=y/n$ with sufficiently large fixed $y$. In
particular, if $p \rightarrow 0$ and $np \rightarrow \infty$, then
$\,\, \hbox{rank}(Q_{n,p})/n \rightarrow 1$.
\end{remark}

\section{Proof of Lemma \ref{lemma:separation}}

If $p$ is at least $(\ln n)^2 / n'$ then $G(m,p)$ will with
probability at least $1-o(1/n^3)$ have no vertices with degree at
most $\ln \ln n$, in which case the lemma is trivially  true.
Therefore we can assume $p \leq (\ln n)^2 / n'$

If $G_{m}$ fails to be well separated for some $m$ between $n'$
and $n$ there must be a first $m_0$ with this property. We are
going to bound the probability that a fixed $m$ is $m_0$.

Case 1: $m_0=n'$. The probability that $G_n'$ fails to be
well-separated is at most $n^2$ times the probability that any
particular pair of vertices $v$ and $w$ are both of small degree and
at distance at most 2 from each other.

The probability that $v$ has sufficiently small degree is at most
\begin{equation} \label{eqn:lowdegree}
\sum_{i=0}^{\ln \ln n} \binom{n'-1}{i} p^i (1-p)^{{n'}-i} \leq
(1+o(1))\sum_{i=0}^{\ln \ln n} (n'p)^i (1-p)^{n'} \leq \frac{(\ln
n)^{\ln \ln n}}{n^{c \delta}},
\end{equation}

and the same holds for $w$ even if we assume $v$ has small degree
(although the degree of $v$ and that of $w$ aren't quite
independent, we can bound the probability $w$ has small degree by
the probability it has at most $\ln \ln n$ neighbors not including
$v$).

Since a given pair of vertices both being of small degree is a
monotone decreasing graph property, the FKG inequality gives that
the probability of an edge being present between $v$ and $w$ is at
most $p$ even after we condition on them both being small degree.
Similarly, the probability of the existence of an $x$ adjacent to
both $v$ and $w$ is at most $np^2$.  Combining these facts, the
probability of two small degree vertices being close is at most
$$ n^2 (\frac{(\ln n)^{\ln \ln n}}{n^{c \delta}})^2
(p+np^2) \leq \frac{(\ln n)^{2 \ln \ln n} \ln^4 n}{n^{2c
\delta-1}}=o(n^{1-2c\delta+\ep}). $$

Case 2: $n'<m_0=m<n$.  We bound the probability that $m$ satisfies
the somewhat weaker condition that $G_{m-1}$ is well separated but
$G_m$ is not.  Let $w$ be the vertex newly added to the graph. There
are only two ways that the addition of $w$ can cause $G_m$ to lose
well-separatedness: either $w$ serves as the link between two low
degree vertices $v_1$ and $v_2$ that were previously unconnected, or
$w$ is itself a low degree vertex of distance at most 2 from a
previous low degree vertex $v_0$.

Applying (\ref{eqn:lowdegree}) twice, the probability that any
particular $v_1$ and $v_2$ both have low degree  and $w$ is
connected to both of them is at most $p^2(\frac{(\ln n)^{\ln \ln
n}}{n^{c \delta}})^2$.

Again by (\ref{eqn:lowdegree}) the probability that $w$ is of low
degree is at most $\frac{(\ln n)^{\ln \ln n}}{n^{c \delta}}$, as
is the probability that any particular choice of candidate for
$v_0$ has low degree.  By the FKG inequality, the probability that
$w$ and our candidate $v_0$ share a common neighbor given they
both have small degree is at most $np^2$.

Since there are at most $n^2$ choices for $v_1$ and $v_2$ and at
most $n$ choices for $v_0$, applying the union bound over all
these choices, we obtain that the probability $G_{m-1}$ is well
connected but $G_m$ is not  is at most

$$ (\frac{(\ln n)^{\ln \ln n}}{n^{c \delta}})^2(2 n^2p^2)=o(n^{-2c\delta+\ep})$$

Applying the union bound over all possible $m$ (there are at most
$n$ values for $m$), we obtain that  the probability of the
existence of such an $m_0$ is $o(n^{1-2c\delta+\ep})$. The proof
is complete.

\section{Proof of Lemma \ref{lemma:expansion}}

In order to prove the edge expansion property we first show that
almost surely all small subgraphs of $G(n,\frac{c \ln n}{n})$ will
be very sparse.

\begin{lemma} \label{lemma:sparseness}
For fixed $c$ the probability that $G(n,\frac{c \ln n}{n})$ has a
subgraph containing at most $\frac{n}{\ln^{3/2} n}$ vertices with
average degree at least 8 is $O(n^{-4})$
\end{lemma}

\begin{proof}
Let $q_j$ be the probability that a subset of size $j$ has at least
$4j$ edges.  By the union bound, this is at most $\binom{n}{j}$
times the probability of a particular subset having at least $4j$
edges, so

\begin{eqnarray*}
q_j &\leq& \binom{n}{j} \binom{j^2 /2}{4j} p^{4j} \\
&\leq& (\frac {ne}{j})^j  (\frac{e j c \ln n}{8 n})^{4j} \\
&\leq& (\frac{c e^5 j^3 \ln ^4 n }{n^3})^{j} .\\
\end{eqnarray*}

For $j<n^{1/4}$ this gives $q_j \leq n^{-2j}$, while for
$j>n^{1/4}$ we have (using our upper bound on $j$) $q_j \leq (\ln
n)^{-j/2} =o(n^{-5})$.  By adding up over all $j$ at least 2, we
can conclude that  the failure probability is $O(n^{-4})$,
completing the proof.
\end{proof}

Armed with this lemma we can now prove Lemma \ref{lemma:expansion}; we do so in
two cases depending on the value of $p$.

Case 1: $p< \frac{12 \ln n}{n}$:

Suppose that $G$ failed to expand properly. If this is the case,
there must be a minimal subset $S_0$ with fewer than $|S_0|$ edges
leaving it.  If any vertex in $S_0$ were adjacent to no other vertex
in $S_0$, it would have a neighbor outside $S_0$ (since $S_0$
contains no isolated vertices), and dropping it would lead to a
smaller non-expanding set, a contradiction. Therefore every vertex
in $S_0$ has a neighbor in $S_0$.  By the well-separatedness
assumption the vertices of degree at most $\ln \ln n$ are
non-adjacent and share no common neighbors. Thus it follows that at
most half the vertices in $S_0$ are of degree at most $\ln \ln n$.
Since at most $|S_0|$ edges leave $S_0$, it follows that there are
$\Omega(|S_0| \ln \ln n )$ edges between vertices of $S_0$. But by
Lemma \ref{lemma:sparseness} the probability an $S_0$ with this many
edges exists is $O(\frac{1}{n^4})$, competing the proof for this
case.

Case 2: $p \geq \frac{12 \ln n}{n}$: We estimate the probability
that there is a non-expanding small set directly by using the
union bound over all sets of size $i<n \ln^{-3/2} n$. The
probability in question can be bounded from above by

\begin{eqnarray*}
\sum_{i=1}^{n \ln^{-3/2} n}
\binom{n}{i}\binom{in}{i-1}(1-p)^{i(n-i)-(i-1)} &\leq&
\sum_{i=1}^{n \ln^{-3/2} n} n^i
(en)^{i-1} e^{-i n p (1+o(1))} \\
&=& \frac{1}{en} \sum_{i=1}^{n \ln^{-3/2} n} (n^2
e^{-np(1+o(1))})^i.
\end{eqnarray*}

The lower bound on $p$ guarantees that the summand is
$O(n^{-(4+o(1))i})$, so the probability for any $p$ in this range is
$o(\frac{1}{n^4})$. Notice that in this case we do not need the
well-separatedness assumption.

\section{Proof of Lemma
\protect\ref{lemma:goodmatrices}}

Let $C_0$ be the event that $G_m$ is good for every $m$ between $n'$ and $n$.
Let $C_1$ be the event that $G_m$ has at most $\frac{1}{p \ln n}$ vertices of
degree less than 2 for every $m$ between $n'$ and $n$, $C_2$ be the event that
$G_m$ has maximum degree at most $10 n p$ for each $m$, and $C_3$ be the event
that $G_m$ is well separated and a small set expander for every $m$ between
$n'$ and $n$. We have
$${\hbox{\bf P}}(\neg C_0) \leq {\hbox{\bf P}}(\neg C_0 \wedge C_1 \wedge C_2 \wedge C_3)+
{\hbox{\bf P}}(\neg C_1) + {\hbox{\bf P}}(\neg C_2) + {\hbox{\bf
P}}(\neg C_3).$$

It suffices to bound each term on the right hand side separately,
and we will do so in reverse order.

Lemmas \ref{lemma:expansion} and \ref{lemma:separation} together
show that ${\hbox{\bf P}}(\neg C_3)=O(n^{1-2c\delta +\ep})$.

${\hbox{\bf P}}(\neg C_2)$ is at most the expected number of
vertices of degree at least $10 n p$ in $G_n$, which is at most

$$n \binom{n}{10 n p} p^{10 n p} \leq n (e/10p)^{10np}
p^{10np} \leq n e^{-10 np} =o(n^{-4}).$$

To bound $\P( \neg C_1)$, we note that he probability that some
$G_m$ contains a set of vertices of degree less than 2 of  size $s=
p^{-1} \ln^{-1} n$ is bounded from above by the probability that at
least $s$ vertices in $G_n$ each have fewer than 2 neighbors amongst
the vertices of $G_{n'}$, which by Markov's inequality is at most
$$\frac{1}{s} n (n'p(1-p)^{n'-2}+(1-p)^{n'-1})=O(s n^2  p e^{-\delta n p}  = o(n^{-1/3}). $$

It follows by the union bound that $\P( \neg C_1)= o(n^{-1/3})$.

It remains to estimate the first term, which we will do by the union
bound over all $m$.  Since property $C_1$ implies $G_m$ has few
vertices of small degree, it suffices to estimate the probability
$G_m$ contains a non-nice set while still satisfying properties
$C_1$, $C_2$, and $C_3$.

Let $p_j$ be the probability that conditions $C_1, C_2$, and $C_3$
hold but some subset of at most $j$ vertices without isolated
vertices is not nice. Symmetry and the
union bound give that $p_j$ is at most $\binom{%
m}{j}$ times the probability that the three conditions hold and some
fixed set $S$ of $j$ vertices is not nice.  We will do this in three
cases depending on the size of $j$.

Case 1: ${\frac{1 }{p \sqrt{\ln n}}}\leq j \leq k.$

Direct computation of the probability that a fixed set of $j$
vertices has either 0 or 1 vertices adjacent to exactly one vertex in the set gives:

\begin{eqnarray*}
p_j &\leq& \binom{m}{j} (
(1-j p(1-p)^{j-1})^{m}+m j p(1-p)^{j-1}(1-j p(1-p)^{j-1})^{m-1})\\
&\leq& (m e p \sqrt{\ln n})^j ((1-j p(1-p)^{j-1})^{m}+m j
p(1-p)^{j-1}(1-j
p(1-p)^{j-1})^{m-1}) \\
&\leq& (m e p \sqrt{\ln n})^j ((1-j p e^{-j p(1+o(1))})^{m}+m j p
e^{-jp(1+o(1))} (1-j p e^{-j p
(1+o(1))})^{m-1}) \\
&\leq& (m e p \sqrt{\ln n})^j (e^{-m j p(1+o(1)) e^{-j p
(1+o(1))}})(1+m j p e^{-j p (1+o(1))}).
\end{eqnarray*}

It follows from our bounds on $j$ and $p$ that $m j p e^{-jp}$ tends to
infinity, so the second half dominates the last term of the above sum and we
have:

\begin{equation*}
p_j \leq (m e p \sqrt{\ln n})^j (e^{-m j p(1+o(1)) e^{-j p
(1+o(1))}})(2 m j p e^{-j p (1+o(1))}).
\end{equation*}

Taking logs and using $\delta n \leq m \leq n$ gives:
\begin{eqnarray*}
\ln(p_j) &\leq& (1+o(1))j (\ln (e n p \sqrt{\ln n}) - \delta n p
e^{-j p
(1+o(1))} - p + {\frac{\ln (2 n j p) }{j}}) \\
&\leq& (1+o(1))j (4\ln(np) - \delta n p e^{-k p (1+o(1))}) \\
&=& (1+o(1)) j (4 \ln (np) - {\frac{\delta n p }{(\ln n)^{{\frac{1 }{2}%
}+o(1)}}}).
\end{eqnarray*}

Since $n p>0.5 \ln n$, taking $n$ large gives that the probability
of failure for any particular $j$ in this range is $o(1/n^4)$, and
adding up over all $j$ and $m$ gives that the probability of a
failure in this range is $o(1/n^2)$.

Case 2: $\frac{10}{2c\delta -1} \leq j \leq \frac{1}{p \sqrt{ \ln
n}}$.

Let $b$ be the number of vertices outside $S$ adjacent to at least
one vertex in $S$, and let $a$ be the number of edges between $S$
and the vertices of $G$ outside $S$.  If $G_m$ is to satisfy the
properties $C_2$ and $C_3$ it must be true that $j \leq a \leq 10 j
n p $.

Next, we note that if $S$ is not nice, then at least $b-1$ of the neighbors of
$S$ must be adjacent to at least two vertices in $S$. This implies that $b \leq
\frac{a+1}{2}$.  It follows that

\begin{equation} \label{eqn:smallj}
p_j \leq \binom{m}{j}(\max_{10j< w \leq 10 j n p} {\hbox{\bf P}}(b
\leq \frac{w+1}{2}|a=w)+ {\hbox{\bf P}}(a \leq 10j) \max_{j \leq w
\leq 10j} {\hbox{\bf P}} (b \leq \frac{w+1}{2}|a=w))
\end{equation}




To bound ${\hbox{\bf P}}(b \leq \frac{w+1}{2}|a=w)$, we fix a set
of  $\frac{w+1}{2}$ vertices and bound  the probability that $w$
consequentially randomly selected vertices were in that set. Using
the union bound over all possible sets of $\frac{w+1}{2}$
vertices, we have

\begin{eqnarray*}
{\hbox{\bf P}}(b \leq \frac{w+1}{2}|a=w) &\leq&
\binom{m-j}{\frac{w+1}{2}} (\frac{w+1}{2(m-j)})^w \\
&\leq& (\frac{2 e (m-j)}{w-1})^{\frac{w+1}{2}} (\frac{w+1}{2(m-j)})^w \\
&\leq& (\frac{4w}{m})^{\frac{w-1}{2}}.
\end{eqnarray*}

This last bound is decreasing in $w$ for the entire range under
consideration (our bounds on $j$ guarantee $w$ is at most $\frac {10
n}{\sqrt{\ln n}}$).  Therefore we can plug in the smallest values of
$w$ in (\ref{eqn:smallj}) to get

\begin{eqnarray*}
p_j &\leq& 3 \sqrt{n} \binom{m}{j}({\hbox{\bf P}}(a<10j) (\frac{e
j}
{m})^{\frac{j}{2}}+(\frac{40 j}{m})^{5j})\\
&\leq& 3 \sqrt {n} (\frac{m e}{j})^j ({\hbox{\bf P}}(a<10j) (\frac{e
j} {m})^{\frac{j}{2}}+(\frac{40 j}{m})^{5j})\\
&\leq& 3 \sqrt{n} ({\hbox{\bf P}}(a<10j) (\frac{e^3
n}{j})^{\frac{j}{2}}+(\frac{130j}{n^{\prime}})^{4j}).
\end{eqnarray*}

$a$ here is the sum of at least $n^{\prime} j(1+o(1))$ independent Bernoulli
variables each with probability of success at least $\frac{c \ln n}{n}$.  We
therefore have
\begin{eqnarray*}
{\hbox{\bf P}}(a \leq 10j) &\leq& \sum_{i=0}^{10j} \binom{n' j}{i}
(\frac{c \ln
n}{n})^i (1- \frac{c \ln n}{n})^{n' j (1+o(1))} \\
&\leq& \sum_{i=0}^{10j} (\frac{c e j \ln n}{i})^i n^{-j c \delta (1+o(1))} \\
&=& n^{-j c \delta (1+o(1))},
\end{eqnarray*}
yielding
\begin{equation*}
p_j \leq 4 \sqrt{n} ((\frac{e^{3/2}}{\sqrt{j} n^{(1+o(1))c\delta-1/2
}})^j+(\frac{130 j}{n^{\prime}})^{4j}).
\end{equation*}

Both terms are decreasing in $j$ in the range under consideration, and plugging
in the lower endpoint of our range gives that $p_j=o(1/n^4)$ for each $j$ and
$m$ in the range.  By the union bound the probability of failure in this range
is $o(1/n^2)$.

Case 3: $2 \leq j \leq \frac{10}{2c \delta -1}$

Let $a$ and $b$ be as in case 2.  We again bound the probability of
failure for any fixed set of vertices by the probability that $b
\leq \frac{a+1}{2}$.

We first note that if condition $C_3$ is to be satisfied then this
inequality must also be satisfied any time $a$ is at most $10j$.
This is because if $a$ is in this range it follows that every vertex
in our set is of degree below $\ln \ln n$, and the
well-separatedness condition then guarantees that each edge leaving
our set must go to a different vertex.

Because of this, we can rewrite equation (\ref{eqn:smallj}) as

\begin{eqnarray*}
p_j &\leq& \binom{m}{j}(\max_{10j< w \leq 10 j n p} {\hbox{\bf P}}(b
\leq \frac{w+1}{2}|a=w) \\
&\leq& 4 \sqrt{n} (\frac{130 j}{n^{\prime}})^{4j},
\end{eqnarray*}

where the second inequality comes from our computations in case 2.
Adding up over all $j$ in this range gives that the probability of
failure in this range is $o(n^{-3})$.

\section{Some Littlewood-Offord-Type Results}

The proof of the remaining two lemmas rely on modifications to the
following lemma of Littlewood and Offord \cite{erdos}:

\begin{lemma} \label{lemma:originalL-O} Let $a_i$ be fixed constants, at least $q$
of which are nonzero. Let $z_1, z_2, \dots z_n$ be random,
independent Bernoulli variables which take on 0 and 1 each with
probability 1/2. Then for any fixed $c$,
\begin{equation*}
{\hbox{\bf P}}(\sum_{i=1}^n a_i z_i = c)=O(q^{-1/2}),
\end{equation*}

where the implied constant is independent of $n$, the $a_i$, and
$c$.
\end{lemma}

The variables we are now considering, however, are not equally
likely to be 0 and 1. Thus we need the following special case of a
more general result in \cite{Hal}.

\begin{lemma}

\label{lemma:modifiedlinearL-O}

Let $a_i$ be fixed constants, at least $q$ of which are nonzero. Let
$z_1, z_2, \dots z_n$ be random, independent Bernoulli variables
which take on 1 with probability $p<{\frac{1 }{2}}$, 0 with
probability $1-p$. Then for any fixed $c$,
\begin{equation*}
{\hbox{\bf P}}(\sum_{i=1}^n a_i z_i = c)=O((qp)^{-1/2}),
\end{equation*}

where the implied constant is absolute.
\end{lemma}

\begin{remark}
The theorem is also true (with near identical proof) if one replaces
the
distribution of the $z_i$ by one with ${\hbox{\bf P}}(z_i=1)={\hbox{\bf P}}%
(z_i=-1)=p, {\hbox{\bf P}}(z_i=0)=1-2p$
\end{remark}

\begin{proof}

Let $r_i$ be Bernoulli random variables taking on 1 with probability
$2p$, and 0 with probability $1-2p$. Let $s_i$ be random variables
taking on 1 and 0 with equal probability, and replace $z_i$ by $r_i
s_i$ (which has the same distribution). By Bayes' inequality we have
\begin{equation*}
{\hbox{\bf P}}(\sum _{i=1}^n (a_i r_i) s_i=c) \leq {\hbox{\bf P}}%
(\sum_{i=1}^n a_i r_i s_i=c | \sum_{i=1}^n r_i \geq qp) + {\hbox{\bf P}}%
(\sum _{i=1}^n r_i<qp)
\end{equation*}

Since ${\hbox{\bf E}}(\sum_{i=1}^n r_i)=2qp$ and $Var(\sum_{i=1}^n
r_i) \leq
2qp$, by Chebyshev's inequality the second term on the right is $%
O((qp)^{-1/2})$. In the first term there are at least $qp$ nonzero $a_i r_i$%
, so the bound follows immediately from the original
Littlewood-Offord lemma.
\end{proof}

The other modified Littlewood-Offord result we need is a similar
modification of the Quadratic Littlewood-Offord lemma in \cite{CTV}:

\begin{lemma}

\label{lemma:modifiedquadraticL-O} Let $a_{ij}$ be fixed constants
such that there are at least $q$ indices $j$ such that for each $j$
there are at least $d$ indices $i$ for which $a_{ij} \neq 0$. Let
$z_1, z_2, \dots z_n$ be as in Lemma \ref{lemma:modifiedlinearL-O}.
Then for any fixed $c$
\begin{equation}  \label{eqn:qref}
{\hbox{\bf P}}(\sum_{i=1}^n \sum_{j=1}^n a_{ij} z_i
z_j=c)=O((qp)^{-1/4}),
\end{equation}

where the implied constant is absolute.
\end{lemma}

\begin{proof}

The proof of Lemma \ref{lemma:modifiedquadraticL-O} relies on the
use of the following application of the Cauchy-Schwartz inequality:

\begin{lemma}

\label{lemma:decoupling} Let $X$ and $Y$ be random variables, and
let $E(X,Y) $ be an event depending on $X$ and $Y$. Let $X^{\prime}$
be an independent copy of $X$. Then
\begin{equation*}
{\hbox{\bf P}}(E(X,Y)) \leq ({\hbox{\bf P}}(E(X,Y) \wedge
E(X^{\prime},Y))^{1/2}
\end{equation*}
\end{lemma}

\begin{proof}

By discretizing, we can reduce to the case when $X$ takes a finite
number of values $x_1, \dots, x_n$. From Bayes' identity we have
\begin{equation*}
{\hbox{\bf P}}(E(X,Y)) = \sum_{i=1}^n {\hbox{\bf P}}(E(x_i,Y)) {\hbox{\bf P}}%
(X=x_i)
\end{equation*}

and
\begin{equation*}
{\hbox{\bf P}}(E(X,Y) \wedge E(X,Y^{\prime})) = \sum_{i=1}^n {\hbox{\bf P}}%
(E(x_i,Y))^2 {\hbox{\bf P}}(X=x_i),
\end{equation*}

and the result follows immediately from Cauchy-Schwartz.
\end{proof}

Without loss of generality we can assume that the $q$ indices $j$
given in the assumptions of our lemma are $1 \leq j \leq q$.

Define $X:=(z_i)_{i>q/2}$, $Y:=(z_i)_{i \leq q/2}$. Let $Q(X,Y)$ be
the quadratic form in (\ref{eqn:qref}), and let $E(X,Y)$ be the
event that that form is 0. By Lemma \ref{lemma:decoupling} we have
\begin{equation*}
{\hbox{\bf P}}(Q(X,Y)=c)^2 \leq {\hbox{\bf
P}}(Q(X,Y)=Q(X^{\prime},Y)=c) \leq {\hbox{\bf
P}}(Q(X,Y)-Q(X^{\prime},Y)=0)
\end{equation*}

Thus it is enough to show the right hand side of this is
$O((qp)^{-1/2})$. To estimate the right hand side, we note that
\begin{equation*}
Q(X,Y)-Q(X^{\prime},Y)=\sum_{j \leq q/2} W_j z_j +
f(X,X^{\prime}),
\end{equation*}

where
\begin{equation*}
W_j=\sum_{i>q/2} a_{ij}(z_i-z_i^{\prime})
\end{equation*}

and $f$ is a quadratic form independent of $Y$. As in Lemma \ref%
{lemma:modifiedlinearL-O}, we next use Bayes' inequality to
condition on the
number of nonzero $W_j$. Let $I_j$ be the indicator variable of the event $%
W_j=0$. We have that
\begin{eqnarray*}
& &{\hbox{\bf P}}(\sum_{j \leq d/2} W_j z_j)=-f(X,X^{\prime})) \\
& & \, \, \, \leq {\hbox{\bf P}}(\sum_{j \leq d/2} W_j z_j)=-f(X,X^{\prime})|%
\sum I_j < q/4) + {\hbox{\bf P}}(\sum I_j \geq q/4)
\end{eqnarray*}

By Lemma \ref{lemma:modifiedlinearL-O} the first term on the right
hand side is $O((qp)^{-1/2})$ for any fixed value of $X$, so it
immediately follows that the same holds true for $X$ random.

For the second term, we note that since $Y$ only involves $q/2$
indices, each index in $Y$ must have at least $q/2$ indices in $X$
with $a_{ij} \neq 0 $. If follows from the remark following Lemma
\ref{lemma:modifiedlinearL-O} that ${\hbox{\bf
E}}(I_j)=O((qp)^{-1/2})$, so ${\hbox{\bf E}}(\sum
I_j)=O(q(qp)^{-1/2})$. By Markov's inequality, the second term is also $%
O((qp)^{-1/2})$, and we are done.
\end{proof}

\section{Proofs of lemmas \protect\ref{lemma:singaug} and \protect\ref%
{lemma:nonsingaug}}

The assumption that the pair $(Q_m,Q_{m+1})$ is normal means that
the rows in $Q_m$ which are entirely 0 have no bearing on the rank
of $Q_{m+1}$.  Thus without loss of generality we can drop those
rows/columns and assume that $A$ has no rows which are all 0, at
which point $A$ will still be singular in Lemma
\ref{lemma:singaug}, but will have become nonsingular in Lemma
\ref{lemma:nonsingaug}.

{\bf \noindent Proof of Lemma \ref{lemma:singaug}.}

If the new column is independent from the columns of $A$, then the
rank increases by two after the augmentation (since $A$ is
singular and the matrices are symmetric). Thus if the rank fails
to increase by two then adding a new column does not increase the
rank.

Assume, without loss of generality, that the rank of $A$ is $D$
and the first $D$ rows $x_1, \dots, x_D$ of $A$ are independent.
Then the last row $x_m$ can be written as a linear combination of
$A$ in a unique way

$$x_m =\sum_{i=1}^D a_i x_i. $$

By throwing away those $a_i$ which are zero, we can assume that
there is some $D' \le D$ such that

$$x_m =\sum_{i=1}^{D'} a_i x_i, $$

\noindent where all $a_i \neq 0$. If $D'+1 < k$, then there is a
vertex $j$ which is adjacent to exactly one vertex from $S=\{1,
\dots, D', m\}$, thanks to the goodness of $Q_m$. But this is a
contradiction as the $j$th coordinates of $x_m$ and $\sum_{i=1}^{D'}
a_i x_i$ do not match (exactly one of them is zero). Thus we can
assume that $D' \ge k-1$.

Now look at the new column $(y_1, \dots, y_m)$. Since the rank
does not increase, we should have

 $$x_m' = \sum_{i=1}^{D'} a_i x_i, $$

 \noindent where $x_i'$ is the extension of $x_i$. This implies

 $$y_m =\sum_{i=1}^{D'} a_i y_i. $$

\noindent Since all $a_i$ are non zero, by  Lemma \ref
{lemma:modifiedlinearL-O} the probability that this happens is
$O((Dp)^{-1/2}) = O((kp)^{-1/2})$, concluding the proof.

{\bf \noindent Proof of Lemma \ref{lemma:nonsingaug}.} Let $A$ be a
good non-singular symmetric matrix of order $m$. Let $A^{\prime}$ be
the $m+1$ be $m+1$ symmetric matrix obtained from $A$ by adding a
new random $(0,1)$
column $u$ of length $m+1$ as the $m+1$st column and its transpose as the $%
m+1$st row.

Let $x_1, \dots, x_{m+1}$ be the coordinates of $u$; $x_{m+1}$ is
the lower-right diagonal entry of $A^{\prime}$ and is zero. The
determinant for $A^{\prime}$ can be expressed as

\begin{equation*}
(\det A) x_{m+1} +\sum_{i=1}^m \sum_{j=1}^m c_{ij} x_i x_j =
\sum_{i=1}^m \sum_{j=1}^m c_{ij} x_i x_j =Q
\end{equation*}

\noindent where $c_{ij}$ is the $ij$ cofactor of $A$. It suffices
to bound the probability
that $Q=0$. We can do this using Lemma \ref%
{lemma:modifiedquadraticL-O} if we can show that many of the
$c_{ij}$ are nonzero.

Since $A$ is now  nonsingular, dropping any of the columns of $A$ will lead to a $%
m \times m-1$ matrix whose rows admit (up to scaling) precisely one
nontrivial linear combination equal to 0. If any of the rows in that combination are dropped, we will be left with an $%
m-1 \times m-1$ nonsingular matrix, i.e. a nonzero cofactor.

As above, that combination cannot involve between 2 and $k$ rows
(since $A$ is good, any set of between 2 and $k$ rows has at least
two columns with exactly one nonzero entry, and even after a
column is removed there will still be one left).  The combination
will involve exactly 1 row only when the column removed
corresponds to the only neighbor of a degree 1 vertex (which
becomes isolated upon the removal of its neighbor).  But by
assumption there are only $\frac{1}{p \ln n}=O(\frac{n}{\ln^2 n})$
possibilities for such a neighbor.

It follows that for each index $j$ except at most $O(\frac{n}{\ln^2 n})$
indices there are at least $k$ indices $i$ for which $c_{ij} \neq
0$, and we can therefore apply Lemma \ref{lemma:modifiedquadraticL-O} with $%
q=k$ to get that ${\hbox{\bf P}}(Q=0)$ and ${\hbox{\bf P}}(Q=-\det
A)$ are both $O((pk)^{-1/4})$, proving Lemma \ref{lemma:nonsingaug}.

\section{Open Problems and Avenues for Further Research}

Theorem \ref{cor:main} gives that $\frac{\ln n}{n}$ is a threshold
for the singularity of $Q(n,p)$ but it would still be of interest to
describe the sources of singularity once $p$ drops below the
threshhold. For instance, Theorem \ref{theo:main} states that for
$p>\ln n / 2n$ almost surely the only cause of singularity is the
presence of isolated vertices. However, once $p$ drops below this
$\ln n / 2n$ this will no longer be the case, as $G(n,p)$ will begin
to acquire pairs of degree one vertices having a common neighbor
(corresponding to pairs of equal rows in $Q(n,p)$).

As noted in remark \ref{remark:smallerp}, it is still the case at
this point that $\,\, \hbox{rank}(Q_{n,p})/n \rightarrow 1$, and
this will continue to occur until $pn=O(1)$.  For $y$ fixed and
$p=y/n$, we have from consideration of isolated vertices and the
bounds in Lemma \ref{lemma:nearfullrank} that

\begin{equation*}
1-O(\ln y/y) \leq(1+o(1)) \,\, \hbox{rank}(Q_{n,p})/n  \leq
1-e^{-y} =1-i(G)/n
\end{equation*}

It seems likely that ${\hbox{\bf E}} (\,\,
\hbox{rank}(Q_{n,y/n}))/n$ tends to some function $g(y)$ as $n
\rightarrow \infty$, and it would be of interest to compute $g$.
Azuma's inequality applied to the vertex exposure process guarantees
the ratio is highly concentrated around this $g$, whatever it may
be.

Let us now consider the case when $p$ is above the threshold $\ln n/n$%
. What is the probability that $Q(n,p)$ is singular ? The current
proof gives bounds which tends to zero rather slowly. For $p >
n^{-\alpha}$ we can prove the singularity probability is $O(
n^{-1/4(1-2\alpha)})$. However, for $p < n^{-1/2}$, we can only
prove $O( (\ln \ln n)^{-1/4})$. While it is certain that these
bounds can be improved by tightening the arguments, it is not clear
how to obtain a significant improvement. For instance, we conjecture that in the case $%
p=\Theta (1)$, the singular probability is exponentially small. Such
bounds are known for non-symmetric random matrices \cite{KKS, TV1,
TV2}, but the proofs do not extend for the symmetric case.

We think that Corollary \ref{cor:main} is relatively easy to extend
to models of random graphs where the edges are independent, but are
included with different probabilities.  However, the assumption of
independence between the edges of $G$ seems crucial. In particular,
the results in this paper do not yet apply to the model of random
regular graphs.

\begin{question}

For what $d$ will the adjacency matrix of the random $d$-regular
graph on $n$ vertices almost surely be nonsingular?
\end{question}

For $d=1$, the matrix is trivially non-singular. For $d=2$, the
graph is union of cycles and the matrix will almost surely be
singular (any cycle of length a multiple of 4 leads to a singular
adjacency matrix).  We conjecture that for $3 \le d \le n-4$, the
matrix is again almost surely nonsingular.


\begin{thebibliography}{9}

\bibitem{CTV} K. Costello, T. Tao and V. Vu,  Random symmetric matrices are
almost surely non-singular, to appear, {\it Duke Math J.}.

\bibitem{erdos} P. Erd\"os, On a lemma of Littlewood and Offord, {\it Bull.
Amer. Math. Soc.} {\bf 51} (1945), 898--902.

\bibitem{Hal} G. Hal\'asz, Estimates for the concentration function
of combinatorial number theory and probability, {\it %
Period. Math. Hungar. } {\bf 8} (1977), no. 3-4, 197-211.

\bibitem{Kom1} J. Koml\'os, On the determinant of $(0,1)$ matrices,  {\it %
Studia Sci. Math. Hungar.} {\bf 2} (1967) 7-22.

\bibitem{KKS} J. Kahn, J. Koml\'os, E. Szemer\'edi, On the
probability a random $\pm$ matrix is singular, {\it J. Amer. Math
Soc.} {\bf 8} (1995) 223-240

\bibitem{TV1} T. Tao and V. Vu, On random $\pm 1$ matrices: Singularity and
Determinant, {\it Random Structures and Algorithms} {\bf 28} (2006)
1-23

\bibitem{TV2} T. Tao and V. Vu, On the singularity probability of random Bernoulli matrices,
to appear, {\it J. Amer. Math. Soc}

\end{thebibliography}
\end{document}